\documentclass{article}
\usepackage{graphicx}

\begin{document}

\title{Cuckoo Search: A Brief Literature Review}
\author{Iztok Fister Jr., Xin-She Yang, Du\v{s}an Fister, Iztok Fister \\
Iztok Fister Jr.\\ Faculty of Electrical Engineering and Computer Science, \\
University of Maribor Slovenia. \\
\and
Xin-She Yang \\ 
School of Science and Technology, \\
Middlesex University, United Kingdom. \\
\and 
Du\v{s}an Fister and Iztok Fister \\ 
Faculty of Electrical Engineering and Computer Science,\\
University of Maribor, Slovenia. 
}

\date{}

\maketitle

\abstract{Cuckoo search (CS) was introduced in 2009, and it has attracted great attention due to its promising efficiency in solving many optimization problems and real-world applications. In the last few years, many papers have been published regarding cuckoo search, and the relevant literature has expanded significantly. This chapter summarizes briefly the majority of the
literature about cuckoo search in peer-reviewed journals and conferences found so far.
These references can be systematically classified into appropriate categories, which can be used as a basis for further research.}

{\bf Citation detail:}
I. Fister Jr., X. S. Yang, D. Fister, I. Fister, Cuckoo search: A brief literature review,
in: {\it Cuckoo Search and Firefly Algorithm: Theory and Applications}, 
Studies in Computational Intelligence, vol. 516, pp. 49-62 (2014).


\section{Introduction}
\label{sec:1}
Since the first introduction of Cuckoo Search (CS) by Xin-She Yang and Suash Deb in 2009~\cite{yang2009cuckoo}, the literature of this algorithm has exploded. Cuckoo search, which drew its inspiration from the brooding parasitism of cuckoo species in Nature, were firstly proposed as a tool
for numerical function optimization and continuous problems. Researchers tested this algorithm on some well-known benchmark functions and compared with PSO and GA, and it was found that cuckoo search achieved better results than the results by PSO and GA. Since then, the original developers of this algorithm and many researchers have also applied this algorithm to engineering optimization, where Cuckoo search also showed promising results.

Nowadays cuckoo search has been applied in almost every area and domain of function optimization, engineering optimization, image processing, scheduling, planning, feature selection, forecasting, and real-world applications. A quick search using Google scholar returned 440 papers, while the original paper by Yang and Deb \cite{yang2009cuckoo} has been cited 223 times
at the time of writing of this chapter. A search using Scirus returned 616 hits with 126 journal papers recorded
up to July 2013. While many papers may be still in press, it is not possible to get hold of all these papers.
Consequently, we will focus on the full papers we can get and thus 114 papers are included in this chapter, which
may be one fraction of the true extent of the literature, but they should be representative and useful.

The aim of this chapter is to provide readers with a brief and yet relatively comprehensive list of literature in the last few years. This helps to gain insight into all the major studies concerning this hot and active optimization algorithm. The structure of this chapter is divided in four different parts. Section 2  presents all the main variants of the  cuckoo search variants, including those studies that have been carried out
in numerical and multi-objective optimization. Hybrids algorithms are also included in this part.
Section 3 focuses on engineering optimization, while Section 4 summarizes all the major applications and
their relevant literature. Then, Section 5 discusses implementation and some
theoretical studies. Finally, Section 6 concludes with some suggestions
for further research topics.

\section{Cuckoo Search: Variants and Hybrids}

\subsection{Variants}

The original cuckoo search was first tested using numerical function optimization benchmarks. Usually, this kind of problems represents a test bed for new developed algorithms. In line with this, standard benchmark function suites~\cite{yang2010engineering} have been developed in order to make comparison between algorithms as fair as possible. For example, some original studies in this area are:
\begin{itemize}
\item Cuckoo search via L\'evy flights \cite{yang2009cuckoo}.
\item An efficient cuckoo search algorithm for numerical function optimization \cite{ong2013efficient}.
\item Multimodal function optimisation \cite{jamil2013multimodal}.
\end{itemize}

Cuckoo search can deal with multimodal problems naturally and efficiently. However, researchers have also attempted to improve its efficiency further
so as to obtained better solutions or comparable results to those in the literature~\cite{eiben2003introduction}, and one such study that is worth mentioning is by Jamil and Zepernick~\cite{jamil2013multimodal}.

\begin{figure}[h]
\centering
\includegraphics[width=1.0\textwidth,height=3in]{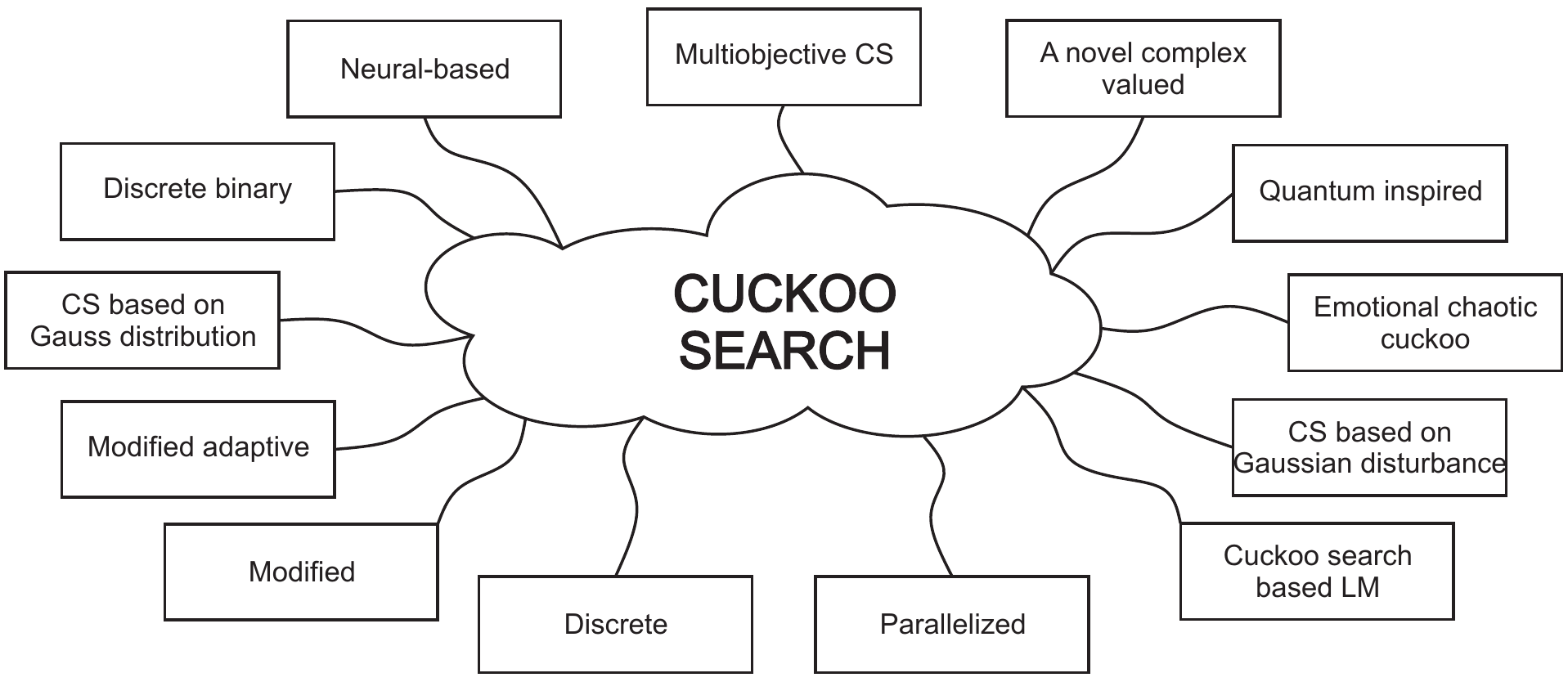}
\caption{Variant of cuckoo search. }
\label{fig:cs}
\end{figure}

Since the first appearance of cuckoo search in 2009, many variants of the cuckoo search algorithm have been developed by
many researchers.  The major variants are summarized in Fig.~\ref{fig:cs} and Table~\ref{tab:mod}.

\begin{table}
\begin{center}
\caption{Variants of Cuckoo Search.}
    \begin{tabular}{ | l | l | l |}
    \hline
    Name & Author & Reference \\ \hline
    Discrete binary CS & Gherboudj et al. & \cite{gherboudj2012solving} \\ \hline
    Discrete CS & Jati and Manurung & \cite{jatidiscrete} \\ \hline
    Discrete CS for TSP & Ouaarab et al. & \cite{ouaarab2013discrete} \\ \hline
    Neural-based CS & Khan and Sahai & \cite{khan2013neural} \\ \hline
    Quantum inspired CS & Layeb & \cite{layeb2011novel} \\ \hline
    Emotional chaotic cuckoo & Lin et al. & \cite{lin2012emotional} \\ \hline
    Cuckoo Search Based LM & Nawi et al. & \cite{Nawi2013}\\ \hline
    Parallelized CS & Subotic et al. & \cite{subotic2012parallelized} \\ \hline
    Modified CS & Tuba et al. & \cite{tuba2011modified} \\ \hline
    Modified CS & Walton et al. & \cite{walton2011modified} \\ \hline
    Modified adaptive CS & Zhang et al. & \cite{zhang2012modified} \\ \hline
    Multiobjective CS & Yang and Deb & \cite{yang2011multiobjective} \\ \hline
    A Novel Complex Valued & Zhou and Zheng & \cite{zhou2013novel} \\ \hline
    CS based on Gauss distribution & Zheng and Zhou  & \cite{zheng2012novel} \\ \hline
    CS based on Gaussian disturbance & Wang et al. & \cite{wang2011cuckoo} \\ \hline
    \end{tabular}
\label{tab:mod}
\end{center}
\end{table}

\subsection{Hybrid Algorithms}
For many continuous optimization problems, cuckoo search can find the desired solutions very efficiently. However, sometimes,
some difficulty may arise, when the appropriate solutions could not be found for some other optimization problems. This is consistent with the so-called No-Free-Lunch theorem~\cite{wolpert1997no}.
To circumvent this theorem, hybridization has been applied to optimization algorithms for solving a given set of problems.
In line with this, cuckoo search has been hybridized with other optimization algorithms, machine learning techniques, heuristics, etc.
Hybridization can take place in almost every component of the cuckoo search. For example,
initialization procedure, evaluation function, moving function and others have all been tried. Some of the hybrid variants
are summarized in Table~\ref{table-hybrid}.
\begin{table}
\begin{center}
\caption{Hybrid cuckoo search. \label{table-hybrid}}
    \begin{tabular}{ | l | l | l |}
    \hline
    Name & Author & Reference \\ \hline
    Hybrid CS/GA & Ghodrati and Lotfi  & \cite{Ghodrati2012, ghodrati2012hybrid} \\ \hline
    Hybrid CS & Li and Yin & \cite{li2013hybrid} \\ \hline
    \end{tabular}
\end{center}
\end{table}

\subsection{Multi-objective Optimization}

Multi-objective optimization consists of more than one objective, and these objectives may be conflicting one another.
Many real-world optimization problems require design solutions according to many criteria.
Single objective optimization searches for a single optimal solution, whilst multi-objective optimization
requires a set of many (potentially infinite), optimal solutions, namely the Pareto front~\cite{robivc2005demo,veldhuizen2000multiobjective}.
Obviously, there are many issues and approaches for multi-objective optimization; however,
two goals in multi-objective optimization are worth noting:
\begin{itemize}
\item to obtain solutions as close to the true Pareto front as possible
\item to generate solutions as diversely as possible in the non-dominated front.
\end{itemize}

Various variants have been developed to extend the standard cuckoo search into multi-objective cuckoo search.
The following list presents some main variants on multi-objective optimization using CS.

\label{sec:MO}
\begin{itemize}
\item Multi-objective CS~\cite{yang2011multiobjective}.
\item Multi-objective scheduling problem~\cite{chandrasekaran2012multi}.
\item Multi-objective cuckoo search algorithm for Jiles-Atherton vector hysteresis parameters estimation~\cite{coelho2013multiobjective}.
\item Pareto archived cuckoo search~\cite{hanoun2012solving}.
\item Hybrid multiobjective optimization using modified cuckoo search algorithm in linear array synthesis \cite{rani2012hybrid}.
\item Multi-objective cuckoo search for water distribution systems~\cite{wang2012multi}.
\end{itemize}

\section{Engineering Optimization}
\label{sec:EO}
\begin{figure}[h]
\centering
\includegraphics[width=1.0\textwidth,height=2.5in]{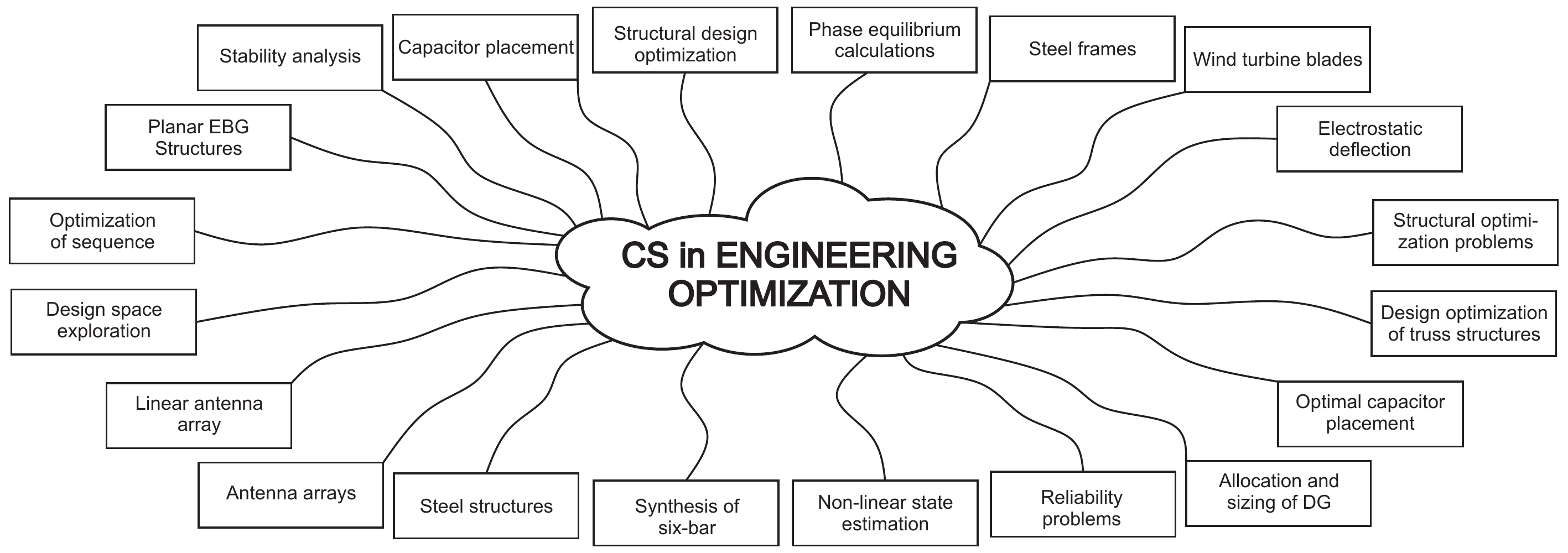}
\caption{Engineering optimization}
\label{fig:cs1}
\end{figure}

Among the diverse applications of cuckoo search, by far the largest fraction of literature
may have focused on the engineering design applications. In fact, cuckoo search and its variants
have become a crucial technology for solving problems in engineering practice as shown in Fig.~\ref{fig:cs1}.
Nowadays, there are applications from almost every engineering domain. Some of these research papers
are summarized in Table~\ref{tab:eng}.

\begin{table}
\begin{center}
\caption{Cuckoo search in engineering optimization}
\label{tab:eng}
    \begin{tabular}{ | l | l | l |}
    \hline
    Problem & Author & Reference \\ \hline
    Engineering optimization & Yang and Deb & \cite{yang2010engineering} \\ \hline
    Capacitor placement & Arcanjo et al. & \cite{arcanjo2012cuckoo} \\ \hline
    Synthesis of six-bar & Bulatovi\v{c} et al. & \cite{bulatovic2013cuckoo} \\ \hline
    Wind turbine blades & Ernst et al. & \cite{ernst2012implementation} \\ \hline
    Design optimization of truss structures & Gandomi et al. & \cite{gandomi2012design} \\ \hline
    Structural optimization problems & Gandomi et al. & \cite{gandomi2013cuckoo} \\ \hline
    Electrostatic deflection  & Goghrehabadi et al. & \cite{goghrehabadi2011hybrid} \\ \hline
    Steel frames & Kaveh and Bakhspoori & \cite{kaveh2011optimum} \\ \hline
    Steel structures & Kaveh et al. & \cite{kaveh2012efficient}\\ \hline
    Antenna arrays & Khodier & \cite{khodier2013optimisation}\\ \hline
    Design space exploration  & Kumar and Chakarverty & \cite{kumar2011design, Kumar2011} \\ \hline
    Optimization of Sequence & Lim et al. & \cite{lim2012cuckoo} \\ \hline
    Planar EBG Structures  & Pain et al. & \cite{pani2013design} \\ \hline
    Stability analysis & Rangasamy and Manickam & \cite{rangasamystability} \\ \hline
    Linear antenna array & Rani and Malek & \cite{rani2011symmetric, rani2012nature} \\ \hline
    Optimal Capacitor Placement & Reddy and Manohar & \cite{reddyoptimal}\\ \hline
    Allocation and sizing of DG & Tan et al. & \cite{tan2012allocation} \\ \hline
    Reliability problems & Valian et al. & \cite{valian2012improved, valian2012cuckoo} \\ \hline
    Non-linear state estimation & Walia and Kapoor & \cite{walia2013particle}\\ \hline
    Phase equilibrium calculations & Bhargava et al. & \cite{Bhargava2013191} \\ \hline
    Structural design optimization & Durgun and Yildiz  & \cite{durgun2012structural}\\ \hline
    \end{tabular}
\end{center}
\end{table}

\section{Applications}
\label{sec:app}
\begin{figure}
\centering
\includegraphics[width=1.0\textwidth,height=3in]{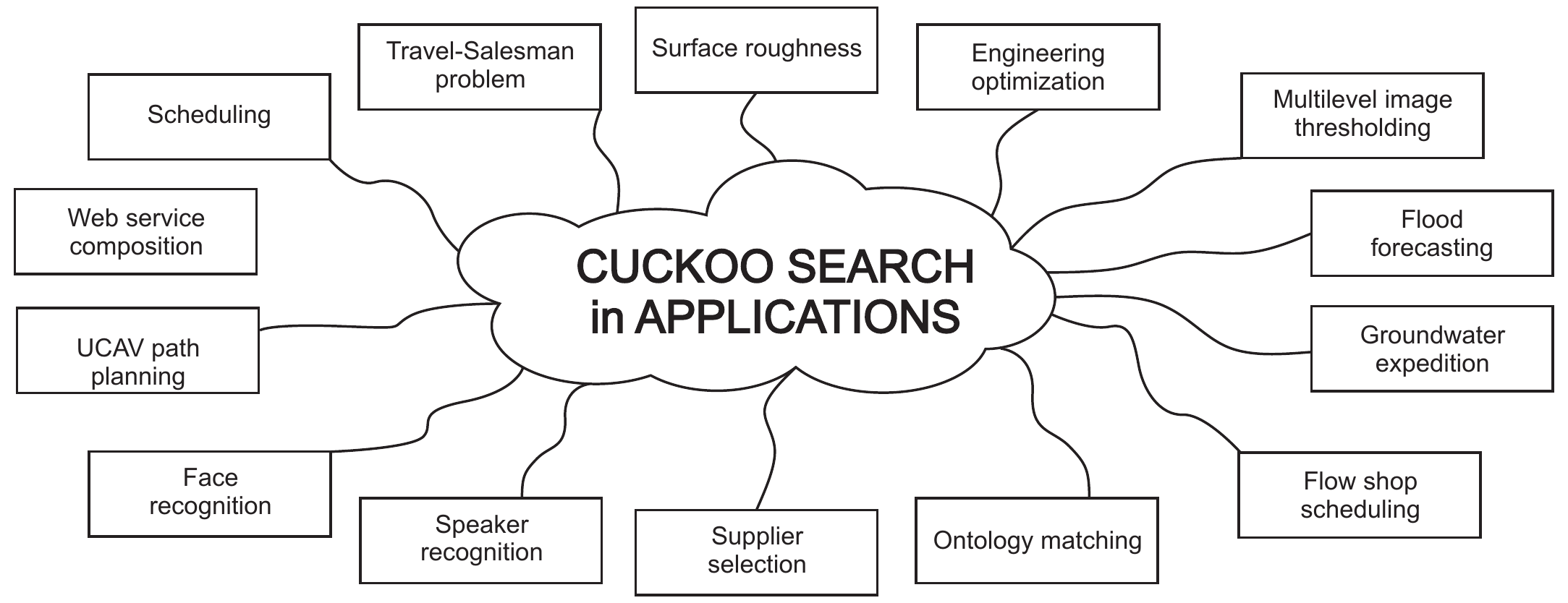}
\caption{Cuckoo search in applications.}
\label{fig:cs2}
\end{figure}

Obviously, engineering optimization is just part of the diverse applications. In fact, cuckoo search
and its variants have been applied into almost every area of sciences, engineering and industry.
Some of the application studies are summarized in Fig.~\ref{fig:cs2} and also in Table~\ref{tab:appls}.

\begin{table}
\begin{center}
\caption{Cuckoo search in applications.}
\label{tab:appls}
    \begin{tabular}{ | l | l | l |}
    \hline
    Application & Author & Reference \\ \hline
    Multilevel image thresholding & Brajevic et al. & \cite{brajevic2012multilevel} \\ \hline
    Flood forecasting & Chaowanawatee \& Heednacram & \cite{chaowanawatee2012implementation} \\ \hline
    Wireless sensor networks & Dhivya \& Sundarambal & \cite{dhivya2011cuckoo} \\ \hline
    Data fusion & Dhivya et al. & \cite{Dhivya2011} \\ \hline
    Cluster in wireless networks & Dhivya et al. & \cite{dhivya2011energy}\\ \hline
    Clustering & Goel et al. & \cite{goel2011cuckoo} \\ \hline
    Groundwater expedition & Gupta et al. & \cite{gupta2013applying} \\ \hline
    Supplier selection & Kanagaraj et al. & \cite{kanagaraj2012supplier} \\ \hline
    Load forecasting & Kavousi-Fard \& Kavousi-Fard & \cite{kavousi2013new} \\ \hline
    Surface Roughness & Madic et al. & \cite{madicapplication} \\ \hline
    Flow shop scheduling & Marichelvam & \cite{marichelvam2012improved} \\ \hline
    Optimal replacement & Mellal et al. & \cite{mellal2012optimal} \\  \hline
    DG allocation in network & Moravej \& Akhlaghi & \cite{moravej2013novel} \\ \hline
    Optimization of Bloom Filter & Natarajan et al. & \cite{natarajan2012enhanced,natarajan2012bloom,natarajan2012comparative} \\ \hline
    BPNN Neural Network & Nawi et al. & \cite{nawi2013new} \\ \hline
    Travelling salesman problem & Ouaarab et al. & \cite{ouaarab2013discrete} \\ \hline
    Web service composition & Pop et al. & \cite{pop2011cuckoo} \\ \hline
    Web service composition & Chifu et al. & \cite{chifu2011bio, chifu2012optimizing} \\ \hline
    Ontology matching & Ritze and Paulheim & \cite{ritze2011towards}\\ \hline
    Speaker recognition & Sood and Kaur & \cite{soodspeaker} \\ \hline
    Automated software testing & Srivastava et al. & \cite{srivastava2012automated, srivastava2012optimal,srivastava2012software} \\ \hline
    Manufacturing optimization & Syberfeldt \& Lidberg & \cite{syberfeldt2012real}\\ \hline
    Face recognition & Tiwari & \cite{tiwari2012face} \\ \hline
    Training neural models & V{\'a}zquez & \cite{vazquez2011training} \\ \hline
    Non-convex economic dispatch & Vo et al. & \cite{vo2013cuckoo} \\ \hline
    UCAV path planning & Wang et al. & \cite{wang2012hybrid,Wang2012} \\ \hline
    Business  optimization & Yang et al. & \cite{yang2012cuckoo} \\ \hline
    Machining parameter selection & Yildiz & \cite{yildiz2013cuckoo} \\ \hline
    Job scheduling in grid & Prakash et al. & \cite{prakashoptimal} \\ \hline
    Quadratic Assignment & Dejam et al. & \cite{dejam2012combining}\\ \hline
    Sheet nesting problem & Elkeran & \cite{elkeran2013new} \\ \hline
    Query optimization & Joshi \& Srivastava & \cite{joshi2013query} \\ \hline
    n-Queens puzzle & Sharma and Keswani & \cite{sharma2013puzzle} \\ \hline
    Computer games & Speed & \cite{speed2010evolving, speed2011artificial} \\ \hline
    \end{tabular}
\end{center}
\end{table}

\section{Theoretical Analysis and Implementation}
As we have seen, the applications of cuckoo search are very diverse. In contrast, the theoretical studies are very limited.
This brief summary may highlight the need for further research in theoretical aspects of cuckoo search.

\subsection{Theory and Algorithm Analysis}
\label{sec:Theory}
It may be difficult to classify a study into a theoretical category or not because the contents
may sometime include both simulations and some analysis of the algorithm. So the following
categorization may not be rigorous. Even so, some theoretical studies about
cuckoo search in the current literature can be summarized, as follows:

\begin{itemize}
\item A conceptual comparison of the cuckoo-search, particle swarm optimization, differential evolution and artificial bee colony algorithms \cite{civicioglu2013conceptual}.
\item Enhancing the performance of cuckoo search algorithm using orthogonal learning method \cite{li2013enhancing}.
\item Starting configuration of cuckoo search algorithm using centroidal Voronoi tessellations \cite{shatnawi2011starting}.
\item Reduced order mesh optimisation using proper orthogonal decomposition and a modified cuckoo search \cite{walton2013reduced, walton2013selected}.
\item Bat algorithm and cuckoo search: a tutorial \cite{yang2013bat}.
\item Metaheuristic algorithms for inverse problems \cite{yang2013metaheuristic,Yang2012a,Yang2012}.
\item Markov model and convergence analysis of cuckoo search \cite{wang2012markov}.
\item Towards the improvement of cuckoo search algorithm \cite{soneji2012towards}.
\end{itemize}

\subsection{Improvements and Other Studies}
\label{sec:Other}
As mentioned earlier, it is not always clear how to classify certain papers. Many research studies concern
the improvements of the standard cuckoo search algorithm. So we loosely put some papers here and thus summarized them
as follows:
\begin{itemize}
\item Tsallis entropy \cite{agrawal2013tsallis}.
\item Improved scatter search using cuckoo search \cite{al2013improved}.
\item Cuckoo search via L{\'e}vy flights for optimization of a physically-based runoff-erosion model \cite{freire2012cuckoo}.
\item Improved differential evolution via cuckoo search operator \cite{musigawan2012improved}.
\item Cuckoo search with the conjugate gradient method \cite{salimi2012extended}.
\item Cuckoo search with PSO \cite{wang2011hybrid}.
\end{itemize}

\subsection{Implementations}
\label{sec:IMP}

Whatever the algorithms may be, proper implementations are very important. Yang provided a
standard demo implementation of cuckoo search\footnote{http://www.mathworks.co.uk/matlabcentral/fileexchange/29809-cuckoo-search-cs-algorithm}.
Important implementations such as object-oriented approach and parallelization have been carried out, as summarized as follows:
\begin{itemize}
\item Object oriented implementation of CS \cite{bacanin2012implementation,bacanin2011object}
\item Parallelization of CS \cite{jovanovicparallelization}.
\end{itemize}

\section{Conclusion}
\label{Conclusion}

In this brief review, a relatively comprehensive bibliography regarding cuckoo search algorithm has been presented.
References have been systematically sorted into proper categories. The rapidly expanding literature implies that
cuckoo search is a very active, hot research area.  There is no doubt that more studies
on cuckoo search will appear in the near future.

From the above review, it is worth pointing out that there are some important issues that need more studies.
One thing is that theoretical analysis should be carried out so that insight can be gained into
various variants of the cuckoo search algorithm. In addition, it may be very useful to carry out parameter
tuning in some efficient variants and see how parameters can affect the behaviour of an algorithm.
Furthermore, applications should focus on large-scale  real-world applications.

\end{document}